\documentclass[12pt,reqno]{amsart}
\usepackage{mathrsfs}
\usepackage{amsgen}
\usepackage{amscd}
\usepackage{amsmath}
\usepackage{latexsym}
\usepackage{amsfonts}
\usepackage{amssymb}
\usepackage{amsthm}
\usepackage{graphicx}
\usepackage{color}
\usepackage{graphicx}
\usepackage{enumerate}
\usepackage{amssymb,amsmath,graphicx,amsfonts,euscript}
\usepackage{mathrsfs}
\usepackage[hyphens]{url}
\usepackage{caption}
\usepackage{subcaption}
\usepackage[margin=3cm, a4paper]{geometry}
\DeclareGraphicsExtensions{.eps}

\newtheorem{theorem}{Theorem}[section]

\newtheorem{problem}[theorem]{Problem}

\theoremstyle{definition}

\newtheorem{assumption}[theorem]{Assumption}

\theoremstyle{remark}

\def\R{\mathbb{R}}

\newcommand{\s}{\mathop{\mathrm{span}}}

\numberwithin{equation}{section}

\allowdisplaybreaks

\begin{document}

\title[]{THE EXPONENTIAL TRAPEZOIDAL METHOD FOR SEMILINEAR INTEGRO-DIFFERENTIAL EQUATIONS}

\author[A. Ostermann]{Alexander Ostermann}
\address{\hspace*{-12pt}A.~Ostermann: Department of Mathematics, Universit\"{a}t Innsbruck, Technikerstrasse 13, 6020 Innsbruck, Austria}
\email{alexander.ostermann@uibk.ac.at}

\author[N. Vaisi]{Nasrin Vaisi}
\address{\hspace*{-12pt}N.~Vaisi: Department of Mathematics, Farhangian University, Sanandaj, Iran}
\email{nasrin\_vaisi@yahoo.com}

\begin{abstract}
The exponential trapezoidal rule is proposed and analyzed for the numerical integration of semilinear integro-differential equations. Although the method is implicit, the numerical solution is easily obtained by standard fixed-point iteration, making its implementation straightforward. Second-order convergence in time is shown in an abstract Hilbert space framework under reasonable assumptions on the problem. Numerical experiments illustrate the proven order of convergence.
\end{abstract}



\maketitle

%%%%%%%%%%%%%%%%%%%%%%%%%%%%%%%%%%%%%%%%%%%%%%%%
\section{Introduction}\label{sec:introduction}
%%%%%%%%%%%%%%%%%%%%%%%%%%%%%%%%%%%%%%%%%%%%%%%%

In this paper we consider the full discretization of an abstract semilinear integro-differential equation of the form
\begin{align}\label{eq1}
u'(t) +\int_{0}^{t}K(t-s)Au(s)\,ds=f(u),  \quad  t \in [0,T],\quad u(0)=u_{0},
\end{align}
where $-A$ is an elliptic differential operator and $K$ is a real-valued positive definite kernel, i.e., for any $T >0$, the kernel $K$ belongs to $L^{1}(0,T)$ and satisfies
\[
\int_{0}^{T} \varphi(t) \int_{0}^{t}K(t-s)\varphi(s)\,ds\,dt\geqslant 0 \quad {\text{for all}} \quad \varphi \in C[0,T].
\]
Equations of the above type and their linear versions are often used to model viscoelastic phenomena and heat conduction in materials with memory. We refer to the monograph~\cite{Pruss} and references therein.

There is an extensive literature on the theoretical and numerical analysis of integro-differential equations \cite{BKO, CTW, LTW, MT, MT2, PTW}. The proposed schemes use finite differences or finite element approximations in space, combined with standard time discretization schemes such as the backward Euler method, the Crank--Nicolson scheme, and other implicit Runge--Kutta or linear multistep methods. The integral term is discretized either by a standard quadrature rule or, in particular for the Riesz kernel, by a convolution quadrature formula \cite{CLP, CP, LST}.

For certain classes of ordinary and partial differential equations, exponential integrators have recently proven to be very efficient. For a survey of these integrators, see~\cite{HO1, HO3, HO2, KT}.  Exponential integrators directly discretize the variation-of-constants formula, which for problem \eqref{eq1} has the form
\begin{equation}\label{eq:voc}
\begin{aligned}
u(t)=S(t)u_{0}+\int^{t}_{0}S(t-\sigma)f(u(\sigma))\,d\sigma,
\end{aligned}
\end{equation}
where $S(t)$ is the solution operator of the linear problem with $f=0$. Exponential integrators can be used to solve this mild form of integro-differential equations. For example, the exponential Euler method applied to~\eqref{eq:voc} is given as
\begin{eqnarray}\label{eq:EE}
U_{m}=S(t_{m})u_{0}+\sum_{j=0}^{m-1}\int^{t_{j+1}}_{t_{j}}S(t_{m}-\sigma)f(U_{j})\,d\sigma, \qquad 1\le m \le M,
\end{eqnarray}
where $U_m$ is an approximation with step size $\tau$ to $u(t)$ at $t=t_m=m\tau$. Note that~\eqref{eq:EE} is an explicit scheme that relies on computing the actions of certain operator functions. The method is efficient if the latter can be done efficiently.

In our previous work \cite{osv}, we proposed explicit exponential Runge--Kutta methods for the time discretization of integro-differential equations. In the linear case, where $f$ is considered only as a function of time, we derived the order conditions for the general order $p$. The resulting exponential quadrature methods were shown to be also convergent of order $p$. In the semilinear case, however, we considered only orders~1 and~2. While the first-order exponential Euler method was simple, the order conditions for the second-order schemes already became involved due to the additional stage required. In this paper, we will consider the exponential trapezoidal integrator as an alternative second-order method for semilinear problems. The method does not require any stages and is easy to implement. Note that, unlike the methods in \cite{osv}, it is implicit. However, since the stiffness of the operator $-A$ is no longer present in the variation-of-constants formula~\eqref{eq:voc}, the resulting system of nonlinear equations can be easily solved by standard methods without any time step restriction due to the stiffness induced by the operator $-A$ in~\eqref{eq1}.

The remainder of this paper is organized as follows. In section 2, we present the abstract framework and some preliminaries. In section 3, we introduce a second-order exponential time integrator for semilinear integro-differential equations along with a spectral Galerkin method for spatial discretization. The error analysis of the proposed integrator is presented in Theorem~\ref{Theorem}, which is the main result of the paper. Finally in section 4, we carry out some numerical experiments and illustrate the theoretical results obtained in the previous sections.

%%%%%%%%%%%%%%%%%%%%%%%%%%%%%%%%%%%%%%%%%%%%%%%%%%%%%%%%%%%%%%%%%%%%%%%
\section{Setting and preliminaries}\label{Setting and preliminaries}
%%%%%%%%%%%%%%%%%%%%%%%%%%%%%%%%%%%%%%%%%%%%%%%%%%%%%%%%%%%%%%%%%%%%%%%

Let $H$ be a real, separable Hilbert space with inner product $(\hspace{0.8mm}\cdot\hspace{0.8mm}{,}\hspace{0.8mm}\cdot\hspace{0.8mm})_H$ and norm $\|v\|_H = \sqrt{(v,v)_H}$. The standard example is $ H=L^2(\mathcal D)$ for a bounded domain $\mathcal D \subset \R^d$, where
$$
(v,w)_{H}=\int_{\mathcal D}v(x)w(x)\,dx,\qquad \|v\|_H=\Big(\int_{\mathcal D} v(x)^{2}\,dx \Big)^{\frac{1}{2}},\qquad v,w\in H.
$$
Furthermore, we denote by $L(H)$ the space of all bounded linear operators on $H$ with the usual operator norm $\| \cdot\|_{L(H)}$. The following assumption will be used.

\begin{assumption}\label{A_operator}
Let $A$ be a densely defined, linear, self-adjoint, positive definite operator on $H$ with compact inverse, and let the kernel $K$ be positive definite.
\end{assumption}

A sufficient condition for $K$ to hold is $k$-monotonicity for some $k\ge 2$. For more details, we refer the reader to \cite{Pruss}, Def.~3.4 and Prop.~3.3. Our prototypical example will be the Riesz kernel
\begin{equation}\label{Rk}
K(t) = \frac{t^{\beta-1}}{\Gamma(\beta)},\qquad t>0,\quad 0<\beta<1.
\end{equation}
However, the framework also includes kernels with less regularity. As an example, we mention a kernel with finite memory, as described in \cite[p.~539]{BGK}.

An important example of $A$ is the negative Laplacian $A =-\Delta$ on a bounded domain $\mathcal{D}\subset\R^{d}$, subject to homogeneous Dirichlet boundary conditions. It is well known that the above assumptions on $A$ imply the existence of a sequence of non-decreasing positive real numbers $\{\lambda_{j}\}_{j=1}^{\infty}$ and an orthonormal basis $\{e_{j}\}_{j=1}^{\infty}$ of $H$ such that
\begin{eqnarray}\label{eq3333}
Ae_{j}=\lambda_{j}e_{j},\quad \lim_{j\rightarrow \infty} \lambda_{j}=\infty.
\end{eqnarray}

For $\nu \in \R$, we consider the domain of $A^{\nu}$, which is a Hilbert space
$$
V=\mathcal{D}(A^{\nu})\quad \text{with norm} \quad \|v\|_{V}=\|A^{\nu}v\|_H.
$$
Our main assumptions on the nonlinearity $f$ are those of \cite{Henry, Pazy}.
\begin{assumption}\label{defination of f}
For some $0\le \nu<1$ and $V=\mathcal{D}(A^{\nu})$, the nonlinearity $f:  V \rightarrow H$ is locally Lipschitz continuous in a neighborhood of the exact solution, i.e.~there exist constants $R>0$ and $L=L(R)$ such that
\begin{equation}\label{eq:lip-f}
\|f(v) - f(w)\|_H \leqslant L\|v - w\|_{V}
\end{equation}
for all $0\le t\le T$ and all $v,w\in V$ satisfying $\|v-u(t)\|_V$, $\|w-u(t)\|_V\leqslant R$.
\end{assumption}

%%%%%%%%%%%%%%%%%%%%%%%%%%%%%%%%%%%%%%%%%%%%%%%%%%%%%%
\subsection{Solution operator.} \label{sec:solop}
A family $\{S(t)\}_{t\geq 0}$ of bounded linear operators on $H$ is called a resolvent family for~(\ref{eq1}) whenever the solution operator $S(t)$ is strongly continuous on $\R_{+}$ and the resolvent equation holds
\[
S(t)u_{0}+\int_{0}^{t}\int_s^t K(\xi-s)\,d\xi\, AS(s)u_{0}\,ds=u_{0},\quad  \text{for all }  u_{0}\in H, \ t\geqslant 0.
\]
If $t\rightarrow u(t)=S(t)u_{0}$ is differentiable for $t > 0$, then $u$ is the unique solution of
\[
u'(t)+\int_{0}^{t}K(t-s)Au(s)\,ds=0,\quad t>0,\quad u(0)=u_{0}.
\]
We refer to the monograph \cite{Pruss} for a comprehensive theory of resolvent families for Volterra equations. Note that the operator family $\{S(t)\}_{t\geq 0}$ does not possess the semigroup property because of the presence  of the memory term in~(\ref{eq1}). Nevertheless, it can be written explicitly using the spectral decomposition (\ref{eq3333}) of $A$ as
\begin{eqnarray}\label{explicit_representation}
S(t)v=\sum_{k=1}^{\infty}s_{k}(t)\,(v,e_{k})_H\,e_{k},
\end{eqnarray}
where the functions $s_{k}(t)$, $k=1,2,\ldots$  are the solutions of the scalar problems
\begin{align}\label{equation}
s'_{k}(t)+\lambda_{k}\int_{0}^{t}K(t-\sigma)s_{k}(\sigma)\,d\sigma=0, \quad s_{k}(0)=1,\qquad  t>0.
\end{align}

Our convergence analysis below will make use of certain smoothing properties of the solution operators.
\begin{assumption}\label{smoothing}
There exist constants $C$ and $1\le \rho\le 2$ such that for any $0\leqslant \alpha \leqslant \frac{1}{\rho}$ the solution operator satisfies the bound
\begin{align}\label{property}
\|A^{\alpha}S(t)\|_{L(H)}\leq Ct^{-\alpha\rho},\quad  t>0.
\end{align}
\end{assumption}

For the Riesz kernel \eqref{Rk}, the smoothing property with $\rho = \beta+1$ is verified in~\cite[Thm~5.5]{MT}. For 3-monotone kernels, estimate~\eqref{property} is verified in~\cite[Prop.~2.2]{Kov14} for
$$
\rho = 1+\frac2{\pi} \sup\left\{\bigl|\arg \widehat K(z)\bigr|\;;\; \text{Re}\, z > 0\right\} \in (1,2),
$$
where $\widehat K$ denotes the Laplace transform of $K$; see also~\cite[Lem.~A.4]{BGK} and \cite[Prop.~3.10]{Pruss}.

%%%%%%%%%%%%%%%%%%%%%%%%%%%%%%%%%%%%%%%%%%%%%%%%%%%%%%%%%%%%%%%%%%%%%%%%
\section{Numerical scheme and main result}\label{Semilinear_problems}
%%%%%%%%%%%%%%%%%%%%%%%%%%%%%%%%%%%%%%%%%%%%%%%%%%%%%%%%%%%%%%%%%%%%%%%%

\subsection{ The numerical method.}
We are now in a position to construct a fully discrete scheme for the numerical solution of problem \eqref{eq1}.
For the spatial discretization we will use a spectral Galerkin method and for the temporal discretization the exponential trapezoidal rule.

For $M\in \mathbb{N}$ we consider the uniform mesh $0=t_{0} < t_{1} < \cdots <t_{M}=T$ on the time interval $\left[0,T\right]$ with time step $h=t_{m+1}-t_{m}$, $m=0,1,\dots,M-1$. Then, by using the variation-of-constants formula, we consider the mild formulation of \eqref{eq1}, viz.
\begin{eqnarray}\label{variation of constant formula1}
u(t_{m})=S(t_{m})u_{0}+\int^{t_{m}}_{0}S(t_{m}-\sigma)f(u(\sigma))\,d\sigma.
\end{eqnarray}
Obviously, this can also be written as
\begin{eqnarray}\label{variation of constant formula}
u(t_{m})=S(t_{m})u_{0}+\sum_{j=0}^{m-1}\int^{t_{j+1}}_{t_{j}}S(t_{m}-\sigma)f(u(\sigma))\,d\sigma.
\end{eqnarray}
Here, the operator $S(t)$ denotes the solution operator of the linear problem (i.e.~for the case $f=0$). We recall that the operator $S(t)$ does not enjoy the semigroup property due to the non-locality of the kernel in \eqref{eq1}.

For the time discretization of \eqref{variation of constant formula}, we employ the \textit{exponential trapezoidal rule}, i.e.
\begin{eqnarray}\label{Trapezoid Integrator3}
U_{m}=S(t_{m})U_{0}+\frac{1}{2}\sum_{j=0}^{m-1}\int^{t_{j+1}}_{t_{j}}S(t_{m}-\sigma)\,d\sigma \lbrace f(U_{j})+f(U_{j+1})\rbrace,
\end{eqnarray}
where $U_m$ ($1\le m \le M$) denotes the numerical approximation to the exact solution $u(t)$ at time $t=t_m$; for notational convenience we set $U_0=u_0$.

For the spatial discretization, we choose $N\in \mathbb N$ and consider the finite dimensional subspace $H_{N} \subseteq H$, given by $H_{N}\equiv \s \lbrace e_{1},e_{2},\cdots,e_{N}\rbrace$, where $\lbrace e_{k}\rbrace_{k=1}^{\infty}$ are the eigenfunctions of $A$, i.e., $Ae_{k}=\lambda_{k}e_{k},~ k\in \mathbb{N}$.
Further we use the projectors $\mathcal{P}_{N}:H\rightarrow H_{N}$ given by
$$
\mathcal{P}_{N}v=\sum_{k=1}^{N}(v,e_{k})e_{k}
$$
for $v \in H$ and the projected operator $A_{N}:H_{N}\rightarrow H_{N}$, $A_{N}=A\mathcal{P}_{N}$
which generates a family of resolvent operators $ \lbrace S_{N}(t)\rbrace_{t\geqslant0}$ in $H_{N}$. It is clear that
\begin{eqnarray}\label{Resolvent_operators}
S_{N}(t)\mathcal{P}_{N}=S(t)\mathcal{P}_{N},
\end{eqnarray}
and also
\begin{eqnarray}\label{Bound}
\Vert A^{-\nu}(I-\mathcal{P}_{N})\Vert_{L(H)}=\sup_{k\geqslant N+1}\lambda_{k}^{-\nu}=\lambda_{N+1}^{-\nu},\quad \nu\geqslant 0.
\end{eqnarray}
A representation of $S_{N}$ is given by
\begin{equation}\label{eq-sn}
S_{N}(t)v=\sum_{k=1}^{N}s_{k}(t)(v,e_{k})e_{k}.
\end{equation}
This motivates us to consider the following fully discrete scheme
\begin{eqnarray}\label{Tra_fully_discrete}
U_{m}^{N}=S_{N}(t_{m})\mathcal{P}_{N}u_{0}+\frac{1}{2}\sum_{j=0}^{m-1}\int^{t_{j+1}}_{t_{j}}S_{N}(t_{m}-\sigma)\,d\sigma\, \mathcal{P}_{N}\Big\{ f(U_{j}^{N})+f(U_{j+1}^{N})\Big\},
\end{eqnarray}
which we propose for the numerical solution of \eqref{variation of constant formula1}.

In order to get a solution in $V$, we assume that the initial data satisfies $u_0\in V$. More regularity, however, improves the spatial convergence result. To elaborate this, we make the following regularity assumption.
\begin{assumption}\label{assumption3b}
Let $g:[0, T] \rightarrow H : t \mapsto g(t) = f(u(t))$ be twice differentiable, let $\nu$ be given by Assumption~\ref{defination of f} and assume that the following conditions hold:
\begin{itemize}\setlength\itemsep{1mm}
\item[(a)] \ $\nu \rho <1$ for $\rho$ given by Assumption~\ref{smoothing};
\item[(b)] \ $u_0=u(0)\in \mathcal D(A^{\nu+\beta})$ for some $\beta \ge 0$;
\item[(c)] \ $A^{\gamma}g \in L^\infty(0,T;H)$ for some $\gamma \ge 0$;
\item[(d)] \ $A^{\eta}g' \in L^\infty(0,T;H)$ for some $0\le \eta \le \nu$;
\item[(e)] \ $A^{-\delta}g''\in L^\infty(0,T;H)$ for some $0\le \delta \le \tfrac1\rho - \nu$.
\end{itemize}
Note that the properties (b)--(e) can also be seen as the definition of the four non-negative parameters $\beta$, $\gamma$, $\delta$, and $\eta$.
\end{assumption}
Under this assumption, we have the following convergence result.
\begin{theorem}\label{Theorem}
For the solution of \eqref{eq1} in the mild form \eqref{variation of constant formula1}, consider the exponential integrator \eqref{Tra_fully_discrete}. If the Assumptions~\ref{A_operator}, \ref{defination of f}, \ref{smoothing}, and~\ref{assumption3b} hold, then there exist constants $h_0>0$ and $C>0$ such that for all step sizes $0<h\le h_0$ and all $N\in \mathbb N$, the global error satisfies for $0<t_m = mh\le T$ and $0\le \alpha<\frac1{\rho}$ the bound
\begin{align*}
\| u(t_{m})-U^{N}_{m}\|_{V}&\leqslant C\Big(t_m^{-\alpha\rho}\lambda_{N+1}^{-\alpha-\beta} + h t_m^{-\nu\rho}\lambda_{N+1}^{-\beta}+\lambda_{N+1}^{\nu-\alpha-\gamma}\\
&\qquad\quad +h^{2-(\nu-\eta)\rho}\sup _{0\leqslant t\leqslant T}\| A^{\eta}g^{\prime}(t)\|_{H}+h^{2}\sup _{0\leqslant t\leqslant T}\| A^{-\delta}g^{\prime\prime}(t)\|_{H} \Big),
\end{align*}
where the constant $C$ depends on $T$, but it is independent of $N$, $m$, and $h$.
\end{theorem}
In particular, if $g'$ is uniformly bounded in $V$, we can choose $\eta=\nu$ and the scheme turns out to be second-order convergent in time.
%%%%%%%%%%%%%%%%%%%%%%%%%%%%%%%%%%
%%%%%%%%%%%%%%%%%%%%%%%%%%%%%%%%
\section{Proof of Theorem~\ref{Theorem}}
First recall that $g(t) = f(u(t))$ and that the operators $A$, $S$, and $P_N$ commute.
By subtracting the numerical solution \eqref{Tra_fully_discrete} from the exact solution \eqref{variation of constant formula} we have
\begin{multline*}
u(t_{m})- U_{m}^{N} =S(t_{m})u_{0}-S_{N}(t_{m})\mathcal{P}_{N}u_{0}+\sum_{j=0}^{m-1}\int^{t_{j+1}}_{t_{j}}\Big\{ S(t_{m}-\sigma)f(u(\sigma))\\
-S_{N}(t_{m}-\sigma)\mathcal{P}_{N} \Big(\frac{1}{2}\big( f(U_{j}^{N})+f(U_{j+1}^{N})\big)\Big)\Big\}\,d\sigma,
\end{multline*}
which, by \eqref{Resolvent_operators}, can be written as
\begin{align*}
u(t_{m})- U_{m}^{N} &=S(t_{m})( I - \mathcal{P}_{N})u_{0}\\
&\quad +\int^{t_{m}}_{0}S(t_{m}-\sigma) \big( g(\sigma)-\mathcal{P}_{N} g(\sigma)\big)\, d\sigma\\
&\quad +\sum_{j=0}^{m-1}\int^{t_{j+1}}_{t_{j}}S(t_{m}-\sigma)\mathcal{P}_{N}
\Big\{ g(\sigma)
- \frac{1}{2}\Big(f(U_{j}^{N})+f(U_{j+1}^{N})\Big)\Big\} \,d\sigma.
\end{align*}
Now taking the norm in $V$, we obtain
\begin{align*}
\Vert u(t_{m})- U_{m}^{N}\Vert_{V} &\leqslant \big\| S(t_{m})( I - \mathcal{P}_{N})u_{0}\big\|_{V}\\
&\quad + \int^{t_{m}}_{0}\Big\| S(t_{m}-\sigma) \big( g(\sigma)-\mathcal{P}_{N} g(\sigma)\big)\Big\|_{V} d\sigma \\
&\quad + \sum_{j=0}^{m-1}\int^{t_{j+1}}_{t_{j}}\Big\| S(t_{m}-\sigma)\mathcal{P}_{N}
\Big\{ g(\sigma)- \frac{1}{2}\Big(f(U_{j}^{N})+f(U_{j+1}^{N})\Big)\Big\}\Big\|_{V} \,d\sigma\\
&=I_{1}+I_{2}+I_{3},
\end{align*}
where $I_{1}, I_{2}$ and $I_{3}$ correspond to the spatial and temporal discretization errors respectively.

First, using (\ref{property}) and the fact (\ref{Bound}) enable us to bound $I_{1}$ as follows
\begin{align*}
I_{1} & = \| A^{\alpha}S(t_{m})A^{-\alpha-\beta}(I-\mathcal{P}_{N})A^{\beta}u_{0}\|_{V} \\
&\leqslant \| A^{\alpha}S(t_{m})\|_{L(H)}\| A^{-\alpha-\beta}(I-\mathcal{P}_{N})A^{\nu+\beta}u_{0}\|_{H} \\
&\leqslant Ct_m^{-\alpha\rho}\lambda_{N+1}^{-\alpha-\beta}\Vert A^{\nu+\beta}u_{0}\Vert_{H}\\
&\leqslant Ct_m^{-\alpha\rho}\lambda_{N+1}^{-\alpha-\beta}.
\end{align*}
To estimate $I_{2},$ we again employ (\ref{property}) and (\ref{Bound}) to obtain
\begin{align*}
I_{2}&= \int_{0}^{t_{m}}\| A^{\alpha}S(t_{m}-\sigma)A^{\nu-\alpha-\gamma}(I-\mathcal{P}_{N})A^{\gamma}g(\sigma)\|_{H} \,d\sigma\\
&\leqslant C\int_{0}^{t_{m}} (t_m-\sigma)^{-\alpha\rho} \|A^{\nu-\alpha-\gamma} (I-\mathcal{P}_{N})\|_{L(H)}\,\| A^{\gamma}g(\sigma)\|_{H} \,d\sigma \\
&\leqslant C\lambda_{N+1}^{\nu-\alpha-\gamma}.
\end{align*}
%%%%%%%%%
It remains to estimate the term
\begin{eqnarray}\label{I33}
I_{3}=\sum_{j=0}^{m-1} \int^{t_{j+1}}_{t_{j}}\Big\|S(t_{m}-\sigma)\mathcal{P}_{N}
\Big\{ g(\sigma)-  \frac{1}{2}\Big(f(U_{j}^{N})+f(U_{j+1}^{N})\Big)\Big\} \Big\|_{V}\,d\sigma.
\end{eqnarray}
We put
\[ g(\sigma)=\frac{1}{2}\Big(g(t_{j})+g(t_{j+1})\Big)+g(\sigma)-\frac{1}{2}\Big(g(t_{j})+g(t_{j+1})\Big),\]
on the right-hand side of~(\ref{I33}) to get
\begin{align*}
I_{3} & \leqslant\frac{1}{2}\sum_{j=0}^{m-1}  \int^{t_{j+1}}_{t_{j}}\big\| S(t_{m}-\sigma)\mathcal{P}_{N}\big(g(t_{j})-f(U_{j}^{N})\big)\big\|_{V}\,d\sigma\\
& \quad +\frac{1}{2}\sum_{j=0}^{m-1}\int^{t_{j+1}}_{t_{j}}\big\| S(t_{m}-\sigma)\mathcal{P}_{N}\big(g(t_{j+1})-f(U_{j+1}^{N})\big)\big\|_{V}\,d\sigma\\
&\quad +  \int^{t_{m}}_{t_{m-1}}\Big\|S(t_{m}-\sigma)\mathcal{P}_{N}\Big[g(\sigma)-\frac{1}{2}\Big(g(t_{m-1})+g(t_{m})\Big)\Big]\Big\|_{V}\,d\sigma\\
&\quad + \sum_{j=0}^{m-2}\int^{t_{j+1}}_{t_{j}}\Big\|S(t_{m}-\sigma)\mathcal{P}_{N}\Big[g(\sigma)-\frac{1}{2}\Big(g(t_{j})+g(t_{j+1})\Big)\Big]\Big\|_{V}\,d\sigma\\
&=I_{3,1}+I_{3,2}+I_{3,3}+I_{3,4}.
\end{align*}
%%%%%%%%%%%%%%%
Next we handle these four terms separately. We first note that
$$
\int_{t_j}^{t_{j+1}} (t_m-\sigma)^{-\nu\rho} d\sigma =
\begin{cases}
Ch(t_m-t_j)^{-\nu\rho},&\quad j\le m-1,\\
h (t_m-t_{j+1})^{-\nu\rho},&\quad j< m-1.
\end{cases}
$$
Using the local Lipschitz continuity of $f$ and~(\ref{property}), we infer that
\begin{align*}\label{I3_1}
I_{3,1}&= \frac{1}{2}\sum_{j=0}^{m-1}\int^{t_{j+1}}_{t_{j}} \Big\| A^{\nu}S(t_{m}-\sigma)\mathcal{P}_{N}\big( f(u(t_{j}))-f(U_{j}^{N}) \big) \Big\|_{H} d\sigma\\
& \leqslant C\sum_{j=0}^{m-1}\int^{t_{j+1}}_{t_{j}} \Vert A^{\nu}S(t_{m}-\sigma)\Vert_{L(H)}\,\Vert u(t_{j})- U_{j}^{N}\Vert_{V} \,d\sigma\\
& \leqslant Ch \sum_{j=1}^{m-1}(t_{m}-t_{j})^{-\nu\rho}\| u(t_{j})- U_{j}^{N}\|_V + Ch t_m^{-\nu\rho}\lambda_{N+1}^{-\beta}.
\end{align*}
In the same way, we get
\begin{align*}
I_{3,2}\leqslant Ch\sum_{j=0}^{m-2}(t_{m}-t_{j+1})^{-\nu\rho}\Vert u(t_{j+1}))-U_{j+1}^{N} \Vert_{V} + Ch (t_{m}-t_{m-1})^{-\nu\rho}\Vert u(t_{m}))-U_{m}^{N}\Vert_{V}.
\end{align*}
%%%%%%%%%%%%%%%%%
In order to bound the term
\begin{align*}
I_{3,3}=  \int^{t_{m}}_{t_{m-1}}\Big\|S(t_{m}-\sigma)\mathcal{P}_{N}\Big[g(\sigma)-\frac{1}{2}\Big(g(t_{m-1})+g(t_{m})\Big)\Big]\Big\|_{V}\,d\sigma,
\end{align*}
we need to expand $g$ in a Taylor series with integral remainder as follows
\begin{align*}
g(\sigma)&=\frac{1}{2}\Big(g(t_{m-1})+g(t_{m})\Big)\\
&\quad +\frac{1}{2}\Big(g'(t_{m-1})(\sigma-t_{m-1})+g'(t_{m})(\sigma-t_{m})\Big)\\
&\quad +\frac{1}{2}\int^{\sigma}_{t_{m-1}}(\sigma-\xi_{1})g^{\prime\prime}(\xi_{1}) \,d\xi_{1}-\frac{1}{2}\int_{\sigma}^{t_{m}}(\sigma-\xi_{2})g^{\prime\prime}(\xi_{2}) \,d\xi_{2}.
\end{align*}
Using this in $I_{3,3}$ we arrive at
\begin{align*}
I_{3,3} &\leqslant \frac{1}{2}\int^{t_{m}}_{t_{m-1}}(\sigma-t_{m-1})\big\|A^{\nu-\eta}S(t_{m}-\sigma)\mathcal{P}_{N}A^\eta g'(t_{m-1})\big\|_{H}\,d\sigma\\
 &\quad + \frac{1}{2}\int^{t_{m}}_{t_{m-1}}(t_{m}-\sigma)\big\|A^{\nu-\eta}S(t_{m}-\sigma)\mathcal{P}_{N}A^\eta g'(t_{m})\big\|_{H}\,d\sigma\\
  & \quad+ \frac{1}{2}\int^{t_{m}}_{t_{m-1}}\Big\|A^{\nu+\delta}S(t_{m}-\sigma)\mathcal{P}_{N}\int^{\sigma}_{t_{m-1}}(\sigma-\xi_{1})A^{-\delta}g^{\prime\prime}(\xi_{1}) \,d\xi_{1}\,d\sigma\Big\|_{H}\\
  &\quad+ \frac{1}{2}\int^{t_{m}}_{t_{m-1}}\Big\|A^{\nu+\delta}S(t_{m}-\sigma)\mathcal{P}_{N}\int_{\sigma}^{t_{m}}(\sigma-\xi_{2})A^{-\delta}g^{\prime\prime}(\xi_{2}) \,d\xi_{2}\,d\sigma\Big\|_{H},
\end{align*}
which then yields
\begin{align*}
I_{3,3} \le Ch^{2-(\nu-\eta)\rho}\sup _{0\leqslant t\leqslant T}\| A^{\eta}g^{\prime}(t)\|_{H}+Ch^{2}\sup _{0\leqslant t\leqslant T}\| A^{-\delta}g^{\prime\prime}(t)\|_{H}.
\end{align*}
%%%%%%%%%%%%%%%%%%%%%%%%%%%%%%%
Finally, for estimating
\begin{align*}
I_{3,4}&= \sum_{j=0}^{m-2}\int^{t_{j+1}}_{t_{j}}\Big\| A^{\nu}S(t_{m}-\sigma)\mathcal{P}_{N}\Big[g(\sigma)-\frac{1}{2}\Big(g(t_{j})+g(t_{j+1})\Big)\Big]\Big\|_{H}\,d\sigma\\
&\le \sum_{j=0}^{m-2}(t_{m}-t_{j+1})^{-(\nu+\delta)\rho}\Big\|\int^{t_{j+1}}_{t_{j}}A^{-\delta}\Big[g(\sigma)-\frac{1}{2}\Big(g(t_{j})+g(t_{j+1})\Big)\Big]\Big\|_{H}\,d\sigma,
\end{align*}
we use the following formula (Peano kernel of the trapezoidal rule):
\begin{align*}
\int_{t_{j}}^{t_{j+1}}\Big[g(\sigma)&-\frac{1}{2}\Big(g(t_{j})+g(t_{j+1})\Big)\Big]\,d \sigma \\
&=\int_{t_{j}}^{t_{j+1}}g(\sigma)\,d \sigma-\frac{h}{2}\Big(g(t_{j})+g(t_{j+1})\Big)\\
&=\frac{1}{2}\int_{t_{j}}^{t_{j+1}}(\xi-t_{j})(\xi-t_{j+1})g^{\prime\prime}(\xi) \,d\xi.
\end{align*}
This immediately leads us to
\begin{align*}
I_{3,4}&\le \frac{1}{2}\sum_{j=0}^{m-2}(t_{m}-t_{j+1})^{-(\nu+\delta)\rho}\int^{t_{j+1}}_{t_{j}}(\xi-t_{j})(t_{j+1}-\xi)\big\|A^{-\delta}g^{\prime\prime}(\xi)\big\|_{H} \,d\xi\\
&\leqslant Ch^{2}\sup _{0\leqslant t\leqslant T}\| A^{-\delta} g^{\prime\prime}(t)\|_{H}\,h\sum_{j=0}^{m-2}(t_{m}-t_{j+1})^{-(\nu+\delta)\rho}\\
&\leqslant Ch^{2}\sup _{0\leqslant t\leqslant T}\| A^{-\delta} g^{\prime\prime}(t)\|_{H}.
\end{align*}
%%%%%%%%%%%%%%%%%%%%%%%%%%%%%%%%%%5
Putting the above estimates together implies (for $h$ sufficiently small)
\begin{multline*}
\| u(t_{m})-U^{N}_{m}\|_{V}\le C\Big( h \sum_{j=1}^{m-1}(t_m-t_{j})^{-\nu\rho}\| u(t_{j})- U^{N}_{j}\|_{V} +t_m^{-\alpha\rho}\lambda_{N+1}^{-\alpha-\beta} + \lambda_{N+1}^{\nu-\alpha-\gamma}\\
+ h t_m^{-\nu\rho}\lambda_{N+1}^{-\beta} + h^{2-(\nu-\eta)\rho}\sup _{0\leqslant t\leqslant T}\| A^{\eta}g^{\prime}(t)\|_{H}+ h^{2}\sup _{0\leqslant t\leqslant T}\| A^{-\delta} g^{\prime\prime}(t)\|_{H} \Big).
\end{multline*}
Applying a discrete Gronwall lemma finally shows the desired bound.

%%%%%%%%%%%%%%%%%%%%%%%%%%%%%%%%%%%%%%%%%%%%%%%%%%%%%%%%%%%%
\section{Implementation and numerical experiments}\label{3.4}
%%%%%%%%%%%%%%%%%%%%%%%%%%%%%%%%%%%%%%%%%%%%%%%%%%%%%%%%%%%%

In this section we present some numerical experiments to illustrate the error bounds obtained in Theorem~\ref{Theorem}. We carry out the experiments in one space dimension, choosing $\mathcal D=(0,1)$ and $A=-\frac{d^2}{dx^2}$, subject to homogeneous Dirichlet boundary conditions. Thus, the eigenvalues and (normalised) eigenfunctions of $A$ are
$$
\lambda_{k}=k^{2}\pi^{2}\qquad  \text{and}  \qquad e_{k}(x)=\sqrt{2}\sin(k\pi x), \qquad k\geqslant 1.
$$

\subsection{Explicit representation of the solution.}
We recall form section~\ref{sec:solop} that
\begin{eqnarray*}
(S_{N}(t)(v))(x)=\sum_{k=1}^{N}2s_{k}(t)\sin(k\pi x)\int_{0}^{1}\sin(k\pi \xi)v(\xi)\,d\xi,
\end{eqnarray*}
for $x\in (0,1)$. The functions ${s}_{k}$ are the solutions of \eqref{equation}. We consider the following two problems.

%%%%%%%%%%%%%%%%%%%%%%%%%%%%%%%%%%
\begin{problem}\label{Riesz}\rm
Consider the Riesz kernel $K(t)=t^{\rho-2}/\Gamma(\rho-1)$ with $1<\rho<2$. Using the Laplace transform in \eqref{equation}, we get
\[
s_{k}(t)=E_{\rho}(-\lambda_{k}t^{\rho}),
\]
where $E_{a}(z)$ is the one-parameter Mittag--Leffler function, defined as $E_{a}(z) = E_{a,1}(z)$, where
$$
E_{a,b}(z)=\sum_{k=0}^{\infty}\frac{z^{k}}{\Gamma(ak+b)},\quad z \in \mathbb{C},\quad a,b>0.
$$
The numerical approximation $U^{N}_{m}$ at time $t_m = m\tau$ can be written as
\begin{align*}
U^{N}_{m}&=\sum_{k=1}^{N} \biggl(E_{\rho}(-\lambda_{k}t^{\rho})(u_{0},e_{k})_H\, e_k\biggr.\\
&\qquad  \biggl.{}+\frac{1}{2}\sum_{j=0}^{m-1}\int_{t_{j}}^{t_{j+1}} E_{\rho}(-\lambda_{k}(t_{m}-\sigma)^{\rho}) \,d\sigma \,
\Big\{  (f(U^{N}_{j}),e_{k})_H+(f(U^{N}_{j+1}),e_{k})_H \Big\} e_{k}\biggr).
\end{align*}
Fixed-point iteration is used to solve this nonlinear problem. Note that the integrals of the Mittag--Leffler function can be computed by a simple quadrature such as the trapezoidal rule. However, the integrals can also be computed exactly as
\begin{align*}
\int_{t_{j}}^{t_{j+1}} E_{\rho}(-\lambda_{k}(t_{m}-\sigma)^{\rho}) \,d\sigma &=\int_{t_{m-j-1}}^{t_{m-j}} E_{\rho}(-\lambda_{k}\sigma^{\rho}) \,d\sigma\\
&= t_{m-j}E_{\rho,2}(-\lambda_{k}t_{m-j}^{\rho})-t_{m-j-1}E_{\rho,2}(-\lambda_{k}t_{m-j-1}^{\rho}),
\end{align*}
which is proved in \cite[Eq.~(1.100)]{Podlubny}. For evaluating the Mittag--Leffler function we use the routine from \cite{PK}.
\end{problem}

%%%%%%%%%%%%%%%%%%%%%%%%%%%%%%%%%%%%%%%%%%%%%%%%%%%%%%%%%%%
\begin{problem}\label{Exponential}\rm
Let $K$ be the smooth kernel
\[
K(t)=e^{-at}  \quad {\text{with}}\quad 0<a\leqslant 2 \quad {\text{for}}\quad t\geqslant 0.
\]
Since $K' = -aK$, it is easy to see that the ordinary integro-differential equation~\eqref{equation} is equivalent to
\[
s_{k}^{\prime\prime}+as_{k}^{\prime}+\lambda_{k}s_{k}=0,\quad s_{k}(0)=1,\quad s_{k}^{\prime}(0)=0.
\]
It has the solution
\[
s_{k}(t)=e^{-\frac{a}{2}t}\Big\{\cos\sqrt{\frac{4\lambda_{k}-a^{2}}{4}}t+\frac{a}{\sqrt{4\lambda_{k}-a^{2}}}\sin\sqrt{\frac{4\lambda_{k}-a^{2}}{4}}t\Big\}.
\]
The integrals of $s_k$ can be computed exactly.
\end{problem}

%%%%%%%%%%%%%%%%%%%%%%%%%%%%%%%%%%%%%%%%%%%%%
\subsection{Numerical experiments.}
In all experiments, we chose the nonlinearity $f(u)= \sin u$, the initial data $u_{0}=4x(1-x)$, $x\in [0,1]$ and $N=100$ frequencies. The problems were integrated with various time step sizes $h$ up to time $T = Mh = 1$ and the errors were calculated in a discrete $L^{2}$-norm using the difference between the numerical solution $U^N_M$ and a reference solution $U_\text{ref}^N$ at time $T=1$:
$$
\text{error} = \left(\Delta x\sum_{j=1}^N \Big( U^N_M(x_j) - U_\text{ref}^N(x_j) \Big)^2\right)^{1/2}, \qquad x_j = j\Delta x,\quad \Delta x = \tfrac1{N+1}.
$$
The reference solution was computed with the second-order explicit exponential integrator from~\cite{osv} using sufficiently small time steps.

\begin{figure}[!ht]
\includegraphics[scale=0.5]{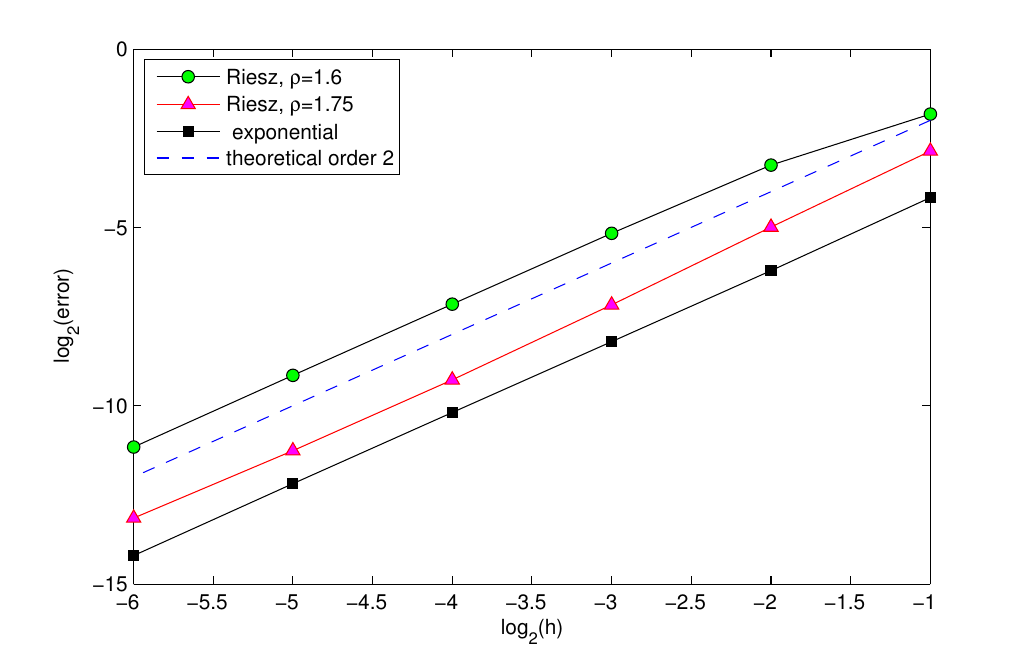}
\caption{\small Temporal rate of convergence of the exponential trapezoidal method for three different problems (see text).\label{fig:err}}
\end{figure}

In the experiments, we considered two different values of $\rho$ for the Riesz kernel and the value $a=2$ for the exponential kernel. Figure~\ref{fig:err} presents a double logarithmic plot of the errors as a function of the time step. The figure confirms the proven theoretical results.

%%%%%%%%%%%%%%%%%%%%%%%%%%%%%%%%%%%%%%%%%%%%%%%%%%%%
\vskip 25pt
\bibliographystyle{model1-num-names}

\end{document}